\documentclass{amsart}
\usepackage{amssymb}
\newcommand{\ri}{{\rm ri\,}}
\newcommand{\cl}{{\rm cl\,}}
\newcommand{\conv}{{\rm conv\,}}
\theoremstyle{plain}
\newtheorem{teo}{Theorem}

\begin{document}
\date{February 26, 2006}
\title{Martingale selection problem
and asset pricing in finite discrete time}
\author{Dmitry B. Rokhlin}
\begin{abstract}
Given a set-valued stochastic process $(V_t)_{t=0}^T$, we say
that the martingale selection problem is solvable if there
exists an adapted sequence of selectors $\xi_t\in V_t$,
admitting an equivalent martingale measure. The aim of this
note is to underline the connection between this problem and
the problems of asset pricing  in general discrete-time market
models with portfolio constraints and transaction costs.
For the case of relatively open convex sets $V_t(\omega)$ we present
effective necessary and sufficient conditions for the solvability of
a suitably generalized martingale selection problem. We show that this
result allows to obtain computationally feasible formulas for the price
bounds of contingent claims. For the case of currency markets we also
give a comment on the first fundamental theorem of asset pricing.

\end{abstract}
\address{Dmitry B. Rokhlin, Faculty of Mechanics and
         Mathematics, Rostov State University, Zorge str.~5,
         Rostov-on-Don, 344090, Russia}
\email{rokhlin@math.rsu.ru}
\keywords{Martingale selection, arbitrage, price bounds, constraints,
transaction costs}
\subjclass{60G42, 91B24}
\maketitle

\section{Introduction}
This paper is motivated by some problems of arbitrage theory.
We deal with discrete-time stochastic securities market models
over general probability spaces.
Recall that in the context of the market model considered
in \cite{DMW90}, the absence of arbitrage opportunities is equivalent
to the existence of an equivalent martingale measure for the discounted
asset price process. Moreover, any arbitrage-free price of a
contingent claim is given by the expectation with respect to some of
these measures. Various generalizations of these results are available.

In spite of their theoretical importance, the purely existence results
of this form are often hard to implement in practice.
Specifically, they are not quite convenient for the calculation of
the price bounds of contingent claims. Another non-trivial issue
concerns the construction of arbitrage strategies in a currency
market with transaction costs.

The present note suggests an approach suitable for these purposes.
In section 2 we present our main tool: the martingale selection
theorem. In section 3 we give two examples, showing that
this result allows to obtain computationally feasible
formulas for the price bounds
of contingent claims in general discrete-time market models with
portfolio constraints and transaction costs.
For the case of currency markets we also make some remarks on
arbitrage opportunities and on the first fundamental theorem of asset
pricing.

\section{Martingale selection theorem}
Consider a probability space $(\Omega,\mathcal F,\mathbf P)$
endowed with the discrete-time filtration $(\mathcal F_t)_{t=0}^T$,
$\mathcal F_0=\{0,\Omega\}$, $\mathcal F_T=\mathcal F$.
In the sequel, all $\sigma$-algebras are assumed to be complete
with respect to $\mathbf P$. We refer to \cite{HP97} for the
definition of measurability (as well as for the definition of a
selector) of a set-valued map.

Suppose $V=(V_t)_{t=0}^T$ is an adapted sequence of
set-valued maps with nonempty relatively open convex values
$V_t(\omega)\subset\mathbb R^d$.
Furthermore, let $(C_t)_{t=0}^{T-1}$ be an adapted sequence of
random convex cones
and let $C_t^\circ$ be the polar of $C_t$.
We say that the {\it $C$-martingale selection problem}
for $(V_t)_{t=0}^T$
is solvable if there exist an adapted stochastic process
$\xi=(\xi_t)_{t=0}^T$ and a probability measure $\mathbf Q$,
equivalent to $\mathbf P$, such that $\xi_t\in V_t$  and
$$\mathbf E_\mathbf Q(\xi_t-\xi_{t-1}|\mathcal F_{t-1})\in C_{t-1}^\circ
\ {\rm a.s.}$$
for all $t\in\{1,\dots, T\}$. Let us call $\xi$ a
$(\mathbf Q,C)$-martingale selector of $V$.
We omit $C$ in all notation if $C_t=\mathbb R^d$.

Given a set $A\subset \mathbb R^d$, denote by
$\cl A$, $\ri A$, $\conv A$ the closure, the
relative interior, and the convex hull of $A$. If $A$
is a cone, then $A^\circ$, $A^*$ are the polar
and the conjugate cones:
$-A^*=A^\circ=\{y:\langle x, y\rangle\le 0,\ \ x\in A\}$.
Here $\langle\cdot,\cdot\rangle$ is the usual scalar product
in $\mathbb R^d$. We also put $A-B=\{x-y:x\in A,\ y\in B\}$.
For a sub-$\sigma$-algebra $\mathcal H\subset\mathcal F$
and a $d$-dimensional $\mathcal F$-measurable random vector $\eta$
denote by $\mathcal K(\eta,\mathcal H;\omega)$ the support of
the regular conditional distribution of $\eta$ with respect to
$\mathcal H$.

Consider an $\mathcal H$-measurable set-valued map $G$
with the closed values $G(\omega)\neq\varnothing$ a.s. Given
a sequence $\{f_i\}_{i=1}^\infty$ of $\mathcal H$-measurable
selectors of $G$ such that the sets $\{f_i(\omega)\}_{i=1}^\infty$ are
dense in $G(\omega)$ a.s. (such a sequence always exists \cite{HP97}),
we put $$\mathcal K(G,\mathcal H;\omega)=\cl\left(\bigcup_{i=1}^\infty
\mathcal K(f_i,\mathcal H;\omega)\right).$$
We refer to \cite{R05}, \cite{R06}, \cite{R06a} for another,
but essentially the
same definition of $\mathcal K(G,\mathcal H)$, which  does not involve
the sequence $\{f_i\}_{i=1}^\infty$ and is expressed
directly in terms of $G$.
If $G(\omega)=\varnothing$ on a set of positive measure, then we put
$\mathcal K(G,\mathcal H)=\varnothing$.

\begin{teo}
The $C$-martingale selection problem for $(V_t)_{t=0}^T$ is solvable
iff the set-valued maps, defined recursively by
$W_T=\cl V_T$;
$$ W_t=\cl(V_t\cap Y_t), \ \
Y_t=\ri(\conv\mathcal K(W_{t+1},\mathcal F_t))-C_t^\circ, \ 0\le t\le T-1$$
have nonempty values a.s. Every $(\mathbf Q,C)$-martingale selector $\xi$
of $V$ take values in $W$ a.s.
\end{teo}

Theorem 1 is an improvement of the main result of \cite{R05},
where the sets $V_t(\omega)$ are assumed to be open and $C_t=\mathbb R^d$.
It is shown in \cite{R06a} that the openness assumption can be dropped.

Sufficiency is the "difficult" part of this theorem. To sketch
the proof, suppose the sets $W_t$ are nonempty and take
some selector $\xi_0\in\ri W_0$. We claim that there exist adapted
sequences $(\xi_t)_{t=0}^T$, $\xi_t\in\ri W_t$;
$(\delta_t)_{t=1}^T$, $\delta_t>0$ such that
$$\mathbf E_\mathbf P(\delta_{t+1}(\xi_{t+1}-\xi_t)|\mathcal F_t)
   \in C_t^\circ;\ \
   \mathbf E_\mathbf P(\delta_{t+1}|\mathcal F_t)=1,\ \
   0\le t\le T-1.$$

These sequences are constructed inductively. Given some $t$,
the induction step is described as follows.
We take some selector $\xi_t\in\ri W_t$ and represent it in the form
$$ \xi_t=\eta_t-\zeta_t,$$
where $\eta_t\in\ri(\conv\mathcal K(W_{t+1},\mathcal F_t))$,
$\zeta_t\in C_t^\circ$ and all elements indexed by $t$
are assumed to be $\mathcal F_t$-measurable. It is crucial
to prove that there exist an element
$\xi_{t+1}\in\ri W_{t+1}$
and a random variable $\delta_{t+1}>0$ such that
$$\eta_t=\mathbf E_\mathbf P(\delta_{t+1}\xi_{t+1}|\mathcal F_t),\ \
  \mathbf E_\mathbf P(\delta_{t+1}|\mathcal F_t)=1.$$
As soon as this is verified (see \cite[Lemma 1]{R06a}), we get
$$ \mathbf E_\mathbf P(\delta_{t+1}(\xi_{t+1}-\xi_t)|\mathcal F_t)=
   \zeta_t\in C_t^\circ.$$

To complete the proof, it remains to introduce the positive
$\mathbf P$-maringale
$$(z_t)_{t=0}^T;\ \ \ z_0=1, \ z_t=\prod_{k=1}^t\delta_k,\ t\ge 1$$
and to check that $\xi$ is a $C$-martingale under the measure $\mathbf Q$
with the density $d\mathbf Q/d\mathbf P=z_N$.

\section{Applications to mathematical finance}
\subsection{Frictionless market with portfolio constraints}
Assume that the discounted prices of $d$ traded assets are described by
a $d$-dimensional adapted stochastic process $(S_t)_{t=0}^T$
and investor's discounted gain is given by
$$G_t^\gamma=\sum_{n=1}^t \langle \gamma_{n-1}, S_n-S_{n-1}\rangle.$$
An adapted admissible portfolio process $\gamma$ is
subject to constraints of the form $\gamma_n\in B_n$,
where $B_n$ are $\mathcal F_n$-measurable random convex cones.
See \cite{PT99}, \cite{EST04}, \cite{R05a} for more information
on this model. The market satisfies the no-arbitrage (NA) condition
if $G_T^\gamma\ge 0$ a.s. implies that $G_T^\gamma=0$ a.s.
for any admissible investment strategy $\gamma$.

Recall, that a contingent claim,
represented by an $\mathcal F_T$-measurable random variable
$f_T$, is called super-hedgeable (resp. sub-hedgeable) at a price
$x\in\mathbb R$ if there exists an admissible portfolio process
$\gamma$ such that $x+G_T^\gamma\ge f_T$ (resp. $x-G_T^\gamma\le f_T$)
a.s. The upper (resp. the lower) price $\overline{\pi}_0$
(resp. $\underline{\pi}_0$) of $f_T$ is the infimum
(resp. the supremum) of all such $x$.

\begin{teo}
Let the cones $B_t$ be polyhedral and assume that NA condition is
satisfied. If the upper and the lower prices of a contingent claim $f_T$
are finite, then they can be computed recursively by
$\overline{\pi}_T=\underline{\pi}_T=f_T;$
$$ \overline{\pi}_t=\sup\{y:(S_t,y)\in\ri(\conv\mathcal
    K((S_{t+1},\overline{\pi}_{t+1}),\mathcal F_t))-B_t^\circ\},\ \
   0\le t\le T-1;$$
$$ \underline{\pi}_t=\inf\{y:(S_t,y)\in\ri(\conv\mathcal
    K((S_{t+1},\underline{\pi}_{t+1}),\mathcal F_t))-B_t^\circ\},\ \
   0\le t\le T-1.$$
\end{teo}

To outline the proof, assume that the contingent
claim $f_T$ is assigned with a price process $(f_t)_{t=0}^T$.
In addition, we allow it to be traded together with $S$ without
additional constraints. It readily follows from known
results that the extended market with the assets $(S_t,f_t)_{t=0}^T$ and
the portfolio constraints $C_t=B_t\times\mathbb R$ is arbitrage-free iff
$(S,f)$ is a $(\mathbf Q,C)$-martingale under some equivalent martingale
measure $\mathbf Q$. The existence of such a process $(f_t)_{t=0}^T$ is
equivalent to the solvability of the $C$-martingale selection problem for
the sequence
$$V_t=\{S_t\}\times\mathbb R,\ t\le T-1;\ \ \ V_T=\{(S_T,f_T)\}.$$

It appears that
$W_t=\{S_t\}\times [\underline\pi_t,\overline\pi_t]$. So,
for any arbitrage-free market extension of the above form we have
$f_t\in [\underline\pi_t,\overline\pi_t]$ a.s. for all $t$.
It can be shown that the bounds of this interval are exactly the upper
and the lower prices of $f_T$ at time $t$.

Theorem 2 is not completely new: for a rather general case of
path-dependent options $f=f(S_0,\dots,S_T)$ after some calculations
it gives the same results as in \cite{CGT01} and \cite{Ru02}.
Certainly, $(\overline\pi_t)_{t=0}^T$ coincides with the minimal
hedging strategy \cite{CGT01}.

\subsection{Currency market with friction}
Our second example concerns Kabanov's model of currency market with
transaction costs \cite{KRS02}, \cite{S04}. Let us briefly describe
this model, literally following \cite{S04}.

Assume that there are $d$ traded currencies.
Their mutual bid and ask prices are specified by an adapted
$d\times d$ matrix process
$(\Pi_t)_{t=0}^T$, $\Pi_t=(\pi_t^{ij})_{1\le i,j\le d}$
such that
$\pi^{ij}>0$, $\pi^{ii}=1$, $\pi^{ij}\le\pi^{ik}\pi^{kj}$.
The solvency cone $K_t$ is generated by the vectors $\{e_i\}_{i=1}^d$
of the standard basis in $\mathbb R^d$ and by the elements
$\pi^{ij}e_i-e_j$. The elements of investor's time-$t$ portfolio
$\theta_t=(\theta_t^i)_{i=1}^d$ represent the amount of each
currency, expressed in physical units. An adapted portfolio process
$\theta=(\theta_t)_{t=0}^T$ is called self-financing if
$\theta_t-\theta_{t-1}\in -K_t$ a.s., $t=0,\dots,T$,
where $\theta_{-1}=0$.

Let $L^0(G,\mathcal H)$ be the set of
$\mathcal H$-measurable elements $\eta$ such that $\eta\in G$ a.s.
Denote by $A_t(\Pi)$ the convex cone in
$L^0(\mathbb R^d,\mathcal F_t)$ formed by the elements $\theta_t$
of all self-financing portfolio processes $\theta$.
According to the definition of \cite{S04}, a bid-ask process
$(\Pi_t)_{t=0}^T$ satisfies the robust no-arbitrage condition
(NA$^r$) if there exists a bid-ask process $(\widetilde\Pi_t)_{t=0}^T$
such that
$$ [1/\widetilde\pi^{ji},\widetilde\pi^{ij}]\subset
\ri [1/\pi^{ji},\pi^{ij}]$$
for all $i,j,t$ and
$ A_T(\widetilde\Pi)\cap L^0(\mathbb R^d,\mathcal F_T)=
  \{0\}.$
An adapted stochastic process $Z=(Z_t)_{t=0}^T$ is a called
strictly consistent price process if $Z$ is a martingale under
$\mathbf P$ and $Z_t\in\ri K_t^*$ a.s., $t=0,\dots,T$.

At first, it is easy to check that the existence of a
strictly consistent price process
is equivalent to the
solvability of the martingale selection problem for the sequence
$(V_t)_{t=0}^T$, $V_t=\ri K_t^*$.
By Theorem 1 this is equivalent to the statement that
the set-valued maps
$W_T=K_T^*$;
$$ W_t=\cl(\ri K_t^*\cap\ri(\conv\mathcal K(W_{t+1},\mathcal F_t))),\ \
0\le t\le T-1$$
have nonempty values a.s. Moreover, it can be shown directly
that the last condition
is equivalent to the following condition of the paper \cite{KRS03}:
if $\sum_{t=0}^T x_t=0$, where $x_t\in L^0(-K_t,\mathcal F_t)$,
then $x_t\in L^0(K_t\cap(-K_t),\mathcal F_t)$, $t=0,\dots,T$.

Together with Theorem 1
this result of \cite{R06} allows
to give a new (but not short) proof for Theorem 1.7 of
\cite{S04} without appealing to the closedness of the cone
$A_T(\Pi)$ of hedgeable claims. Very recently some results in
this direction were also announced by Mikl\'{o}s R\'{a}sonyi.

Furthermore,
if $W_n=\varnothing$ for some $n$, then there exists an element
$x_n\in L^0(-K_n,\mathcal F_n)$ that properly separates the
cones $K_n^*$ and $\conv\mathcal K(W_{t+1},\mathcal F_t)$
on a set $A_n\in \mathcal F_n$ of positive measure.
It is proved in \cite{R06} that $x_n$ admits
a representation
$$x_n=-\sum_{t=n+1}^T x_t,\ \ x_t\in L^0(-K_t,\mathcal F_t)$$
and that there exists a number $m\ge n$ such that
$x_m\not\in K_m\cap(-K_m)$ on a set
of positive measure.
Following \cite{KRS03}, it is straightforward to obtain
an arbitrage portfolio of the form
$$ \theta_{-1}=0,\ \ \theta_t=\theta_{t-1}+x_t+\varepsilon_t,\ \
   0\le t\le T,$$
$\varepsilon_m\in L^0(\mathbb R_+^d,\mathcal F_m)$, $\varepsilon_t=0$,
$t\neq m$ for any bid-ask process $\widetilde{\Pi}$, introduced
above.

The next result describes the set of arbitrage-free endowments
for a contingent claim $\zeta_T\in L^0(\mathbb R^d,\mathcal F_T)$.
\begin{teo}
Suppose NA$^r$ condition is satisfied and $\zeta_0\in\mathbb R^d$ is an
initial endowment. Then there exists a strictly
consistent price process $Z$ such that
$\langle\zeta_0,Z_0\rangle=\mathbf E_\mathbf P\langle\zeta_T,Z_T\rangle$
iff the martingale selection problem for the set-valued stochastic
sequence
\begin{eqnarray*}
V_0 &=& \{(x,y):x\in\ri K_0^*,\ y=\langle \zeta_0,x\rangle\};\ \
V_t =\ri K_t^*\times\mathbb R,\ 1\le t\le T-1;\\
V_T &=& \{(x,y):x\in\ri K_T^*,\ y=\langle \zeta_T,x\rangle\}.
\end{eqnarray*}
is solvable.
The set $\mathcal I(\zeta_T)$ of all such $\zeta_0$ coincides with
$\ri W_0$, where $W_0$ is defined in Theorem 1.
\end{teo}

To gain a better understanding of this result, imagine a frictionless
market $\mathcal M(Y)$ with $d$ traded assets $(Y^i)_{i=1}^d$,
whose exchange rates satisfy the conditions
$$Y_t^j/Y_t^i\in\ri [1/\pi^{ji},\pi^{ij}]$$
or, equivalently, $Y_t\in\ri K_t^*$.
The martingale selection problem for the sequence $(V_t)_{t=0}^T$,
defined in Theorem 3, is solvable iff
$f_0=\langle\zeta_0,Y_0\rangle$ is an arbitrage-free
price of $f_T=\langle\zeta_T,Y_T\rangle$ in the market $\mathcal M(Y)$
for at least one $Y\in\ri K^*$
(see e.g. \cite{S03} for the formal definition of an arbitrage-free price).

Furthermore, suppose $(\xi,\mathbf Q)$ is a solution of this problem, i.e.
$\xi=(Y,f)\in V$ is a martingale under
some probability measure $\mathbf Q$, equivalent to $\mathbf P$.
Denote by $(z_t)_{t=0}^T$ the density process of $\mathbf Q$ with respect
to $\mathbf P$:
$$ z_T=d\mathbf Q/d\mathbf P,\ \ \
   z_{t-1}=\mathbf E_\mathbf P(z_t|\mathcal F_{t-1}),\ \ t\le T.$$
Then $(Z_t,g_t)_{t=0}^T=(z_t Y_t,z_t f_t)_{t=0}^T$ is a martingale
under $\mathbf P$, $Z$ is a strictly consistent price process and
$$ \langle\zeta_0,Z_0\rangle=\langle\zeta_0,Y_0\rangle=f_0=
   \mathbf E_\mathbf Q f_T=
   \mathbf E_\mathbf Q \langle\zeta_T,Y_T\rangle=
   \mathbf E_\mathbf P \langle\zeta_T,Z_T\rangle.$$

At last, we mention that the complement of $\mathcal I(\zeta_T)$
consists of two disjoint convex sets. The closure of one of these sets
contains exactly those endowments, which are needed to superreplicate
$\zeta_T$.

\end{document}